\documentclass{article}
\usepackage{graphicx}
\usepackage{amsmath}
\usepackage{amsfonts}
\usepackage{amssymb}

\begin{document}

\begin{center}
\smallskip

\textbf{Signed degree sets in signed bipartite graphs\bigskip\bigskip}
\end{center}

\textbf{S. Pirzada}$^{1}$,\textbf{\ T. A. Naikoo}$^{2}\medskip$\textbf{ and F.
A. Dar}$^{3}$

Department of Mathematics, University of Kashmir, Srinagar-190006, India

$^{1}$E{mail: sdpirzada@yahoo.co.in}\newline \ \ \ \ \ $^{2}${Email:
tariqnaikoo@rediffmail.com}\newline { \ \ \ \ \ }$^{3}${Email:
sfarooqdar@yahoo.co.in}\newline \medskip

\bigskip

\begin{center}
\textbf{Abstract\bigskip}
\end{center}

A signed bipartite graph G(U, V) is a bipartite graph in which each edge is
assigned a positive or a negative sign. The signed degree of a vertex x in
G(U, V) is the number of positive edges incident with x less the number of
negative edges incident with x. The set S of distinct signed degrees of the
vertices of G(U, V) is called its signed degree set. In this paper, we prove
that every set of integers is the signed degree set of some connected signed
bipartite graph. \medskip\medskip

\begin{center}
\bigskip\textbf{1. Introduction}
\end{center}

A signed graph is a graph in which each edge is assigned a positive or a
negative sign. The concept of signed graphs is given by Harary [3]. Let G be a
signed graph with vertex set V = \{v$_{1}$, v$_{2}$,..., v$_{n}$\}. The signed
degree of a vertex v$_{i}$ in G is denoted by sdeg(v$_{i}$) (or simply by
d$_{i}$) and is defined as d$_{i}$ = d$_{i}^{+}$-d$_{i}^{-}$, where 1 $\leq
$\ i $\leq$\ n and d$_{i}^{+}\ $(d$_{i}^{-}$) is the number of positive
(negative) edges incident with v$_{i}.$ A signed degree sequence
$\sigma=[d_{1},d_{2},...,d_{n}]\ $of a signed graph G is formed by listing the
vertex signed degrees in non-increasing order. An integral sequence is
s-graphical if it is the signed degree sequence of a signed graph. Also, a
non-zero sequence $\sigma=[d_{1},d_{2},...,d_{n}]\ $is a standard sequence if
$\sigma$ is non-increasing ,$\ \sum_{i=1}^{n}d_{i}$ is even,$\ d_{1}>0$, each
$|d_{i}|<n$, and $|d_{1}|\geq|d_{n}|$.\newline \smallskip

The following result [1], gives a necessary and sufficient condition for an
integral sequence to be s-graphical, which is similar to Hakimi's result for
degree sequences in graphs [2].\newline \vspace{0.3cm}

\noindent\textbf{Theorem 1.1. }A standard integral sequence $\sigma
=[d_{1},d_{2},...,d_{n}$] is s-graphical if and only if \newline
$\ \ \ \ \ \ \ \ \ \ \ \ \sigma^{/}=[d_{2}-1,d_{3}-1,...,d_{d_{1}%
+s+1}-1,d_{d_{1}+s+2},...,d_{n-s},d_{n-s+1}+1,...,d_{n}+1,]$\newline is
s-graphical for some s, $0\leq s\leq\frac{n-1-d_{1`}}{2}.\smallskip
\medskip\medskip$

The next characterization for signed degree sequences in signed graphs is due
to Yan et al. [8]. \medskip\medskip\newline \vspace{0.3cm}\noindent
\textbf{Theorem 1.2.} A standard integral sequence $\sigma=[d_{1}%
,d_{2},...,d_{n}$] is s-graphical if and only if \newline
$\ \ \ \ \ \ \ \ \ \ \sigma_{m}^{/}=[d_{2}-1,d_{3}-1,...,d_{d_{1}%
+m+1}-1,d_{d_{1}+m+2},...,d_{n-m},d_{n-m+1}+1,...,d_{n}+1,]$\newline is
s-graphical, where m is the maximum non-negative integer such that
$d_{d_{1}+m+1}>d_{n-m+1}$\newline \medskip

\bigskip The set of distinct signed degrees of the vertices of a signed graph
is called its signed degree set.

In [4], Kapoor et al. proved that every set of positive integers is the degree
set of some connected graph and determined the smallest order for such a
graph.\ Pirzada et al. [7] proved that every set of positive (negative)
integers is the signed degree set of some connected signed graph and
determined the smallest possible order for such a signed graph.\medskip

A graph G is called bipartite if its vertex set can be partitioned into two
nonempty disjoint subsets U and V such that each edge in G joins a vertex in U
with a vartex in V, and is denoted by G(U, V). Let G(U, V) be a bipartite
graph with U = \{u$_{1}$, u$_{2}$,..., u$_{p}$\} \medskip and V = \{v$_{1}$,
v$_{2}$,..., v$_{q}$\}. Then, degree of u$_{i}$ (v$_{j}$) is the number of
edges of G(U, V) incident with u$_{i}$ (v$_{j}$) and is denoted by deg u$_{i}$
or simply by d$_{i}$ (deg v$_{j}$ or simply by e$_{j}$). Then, the sequences
$[d_{1},d_{2},...,d_{p}]$ and $[e_{1},e_{2},...,e_{q}]$\ are called the degree
sequences of G(U, V). The set of distinct degrees of the vertices of a
bipartite graph G(U, V) is called its degree set.\bigskip

The criteria for bipartite graphical sequences is given by Gale-Ryser theorem

\begin{center}
$\sum_{i=1}^{p}d_{i}=\sum_{j=1}^{q}e_{j}$, \ $\sum_{i=1}^{k}d_{i}\leq
\sum_{j=1}^{q}\min\{k,\ e_{j}\}$, for $1\leq k\leq p.\medskip$
\end{center}

Pirzada et al. [6] proved that every set of non-negative integers is a degree
set of some bipartite graph. \medskip

A signed bipartite graph is a bipartite graph in which each edge is assigned a
positive or a negative sign. Let G(U, V) be a signed bipartite graph with U =
\{u$_{1}$, u$_{2}$,..., u$_{p}$\} \medskip and V = \{v$_{1}$, v$_{2}$,...,
v$_{q}$\}. Then,

\begin{center}
signed degree of u$_{i}$ is sdeg(u$_{i}$) = d$_{i}$ = d$_{i}^{+}$-d$_{i}^{-}$,
\end{center}

where 1 $\leq$\ i $\leq$\ p and d$_{i}^{+}\ $(d$_{i}^{-}$) is the number of
positive (negative) edges incident with u$_{i}$, and

\begin{center}
signed degree of v$_{j}$ is sdeg(v$_{j}$) = e$_{j}$ = e$_{j}^{+}$-e$_{j}^{-}$,
\end{center}

where 1 $\leq$\ j $\leq$\ q and e$_{j}^{+}\ $(e$_{j}^{-}$) is the number of
positive (negative) edges incident with v$_{j}$.

Clearly, $\left|  d_{i}\right|  \leq q$ and $\left|  e_{j}\right|  \leq p$.
Then, the sequences $\alpha=[d_{1},d_{2},...,d_{p}]$ and $\beta=[e_{1}%
,e_{2},...,e_{q}]$\ are called the signed degree sequences of the signed
bipartite graph G(U, V). Two sequences $\alpha=[d_{1},d_{2},...,d_{p}]$ and
$\beta=[e_{1},e_{2},...,e_{q}]$\ are said to be standard sequences if \ 

Two sequences of integers $\alpha=[d_{1},d_{2},...,d_{p}]$ and $\beta
=[e_{1},e_{2},...,e_{q}]$\ are said to be standard sequences if

\ \ \ \ \ (i) $\alpha$\ is non-zero, (ii) $\alpha$\ is non-increasing and
$\left|  d_{1}\right|  \geq\left|  d_{p}\right|  $, for we may always replace
$\alpha$ and $\beta$ by $-\alpha$ and $-\beta$ if necessary, (iii) $\sum
_{i=1}^{p}d_{i}=\sum_{j=1}^{q}e_{j},$ (iv) $d_{1}>0,$ (v) each $\left|
d_{i}\right|  \leq q,$ each $\left|  e_{j}\right|  \leq p,$ and each $\left|
e_{j}\right|  \leq\left|  d_{1}\right|  $.\ \ \ \ \ \ \ \ \ \ \ \ \ \ \ \ \ \ \ \ \ \ \ \ \ \ \ \ \ \ \ \ \ \ \ \ \ \ \ \ \ \ \ \ \ \ \ \ \ \ \ 

or \ (i) $\beta$\ is non-zero, (ii) $\beta$\ is non-increasing and $\left|
e_{1}\right|  \geq\left|  e_{q}\right|  $, for we may always replace $\alpha$
and $\beta$ by $-\alpha$ and $-\beta$ if necessary, (iii) $\sum_{i=1}^{p}%
d_{i}=\sum_{j=1}^{q}e_{j},$ (iv) $e_{1}>0,$ (v) each $\left|  d_{i}\right|
\leq q,$ each $\left|  e_{j}\right|  \leq p,$ and each $\left|  d_{i}\right|
\leq\left|  e_{1}\right|  $.\ \ 

Also, a signed bipartite graph G(U, V) is said to be connected if each vertex
u $\in$\ U is connected to every vertex v $\in$\ V.\ \medskip

The next result, due to Pirzada and Naikoo [5], is a necessary and sufficient
condition for a pair of integral sequences to be the signed degree sequences
of some signed bipartite graph.\bigskip

\textbf{Theorem 1.3. }Let $\alpha=[d_{1},d_{2},...,d_{p}]$ and $\beta
=[e_{1},e_{2},...,e_{q}]$\ be standard sequences. Then, $\alpha$\ and $\beta
$\ are the signed degree sequences of a signed bipartite graph if and only if
there exist integers r and s with d$_{1}$ = r-s and 0 $\leq$\ s $\leq
\frac{q-d_{1}}{2}$\ such that $\alpha^{/}$\ and $\beta^{/}$\ are the signed
degree sequences of a signed bipartite graph, where $\alpha^{/}$ is obtained
from $\alpha$\ by deleting d$_{1}$\ and $\beta^{/}$\ is obtained from $\beta
$\ by reducing r greatest entries of $\beta$\ by 1 each and adding s least
entries of $\beta$\ by 1 each.\ \medskip

For any two sets X and Y, we denote by X$\bigoplus$Y to mean that each vertex
of X is joined to every vertex of Y by a positive edge.\ \ \ \ \bigskip

\begin{center}
\textbf{2. Main Results\medskip\medskip}
\end{center}

The set S of distinct signed degrees of the vertices of a signed bipartite
graph G(U, V) is called its signed degree set. \medskip

The following result shows that every set of positive integers is a signed
degree set of some connected signed bipartite graph.\medskip

\textbf{Theorem 2.1.} Let d$_{1}$, d$_{2}$,..., d$_{n}$ be positive integers.
Then, there exists a connected signed bipartite graph with signed degree set S
= $\left\{  d_{1},\sum_{i=1}^{2}d_{i},...,\sum_{i=1}^{n}d_{i}\right\}  .$

\textbf{Proof.} If n = 1, then a signed bipartite graph G(U, V) with $\left|
U\right|  $ = $\left|  V\right|  $ = d$_{1}$ and U$\bigoplus$V has signed
degree set S = \{d$_{1}$\}.

For n $\geq$\ 2, construct a signed bipartite graph G(U, V) as follows.

Let $\ \ U=X_{1}\cup X_{2}\cup X_{2}^{/}\cup...\cup X_{n}\cup X_{n}%
^{/},\smallskip$

\ \ \ \ \ $\ \ V=Y_{1}\cup Y_{2}\cup Y_{2}^{/}\cup...\cup Y_{n}\cup Y_{n}^{/},$

with $X_{i}\cap X_{j}=\phi,\ X_{i}\cap X_{j}^{/}=\phi,\ X_{i}^{/}\cap
X_{j}^{/}=\phi,Y_{i}\cap Y_{j}=\phi,\ Y_{i}\cap Y_{j}^{/}=\phi,\ Y_{i}^{/}\cap
Y_{j}^{/}=\phi\ (i\neq j),$

$\left|  X_{i}\right|  =\left|  Y_{i}\right|  =d_{i}$ for all i, 1 $\leq$\ i
$\leq$\ n, $\left|  X_{i}^{/}\right|  =\left|  Y_{i}^{/}\right|  =d_{1}%
+d_{2}+...+d_{i-1}$ for all i, 2 $\leq$\ i $\leq$\ n. Let (i) X$_{i}\bigoplus
$Y$_{j}$ whenever i $\geq$\ j, (ii) X$_{i}^{/}\bigoplus$Y$_{i}$ for all i, 2
$\leq$\ i $\leq$\ n, and (iii) X$_{i}^{/}\bigoplus$Y$_{i}^{/}$ for all i, 2
$\leq$\ i $\leq$\ n. Then, the signed degrees of the vertices of G(U, V) are
as follows.

For 1 $\leq$\ i $\leq$\ n

\ \ \ \ \ sdeg(x$_{i}$) = $\sum_{j=1}^{i}\left|  Y_{j}\right|  =\sum_{j=1}%
^{i}d_{j}=d_{1}+d_{2}+...+d_{i},$ for all x$_{i}$\ $\in$\ X$_{i}$,

for 2 $\leq$\ i $\leq$\ n

\ \ \ \ \ sdeg(x$_{i}^{/}$) = $\left|  Y_{i}\right|  +\left|  Y_{i}%
^{/}\right|  =d_{i}+d_{1}+d_{2}+...+d_{i-1}\smallskip$

$\ \ \ \ \ \ \ \ \ \ \ \ \ \ \ \ \ =d_{1}+d_{2}+...+d_{i},$ for all x$_{i}%
^{/}$\ $\in$\ X$_{i}^{/}$, for 1 $\leq$\ i $\leq$\ n,

\ \ \ \ \ sdeg(y$_{i}$) = $\sum_{j=i}^{n}\left|  X_{j}\right|  +\left|
X_{i}^{/}\right|  =\sum_{j=i}^{n}d_{j}+d_{1}+d_{2}+...+d_{i-1}\smallskip$

$\ \ \ \ \ \ \ \ \ \ \ \ \ \ \ \ \ =d_{i}+d_{i+1}+...+d_{n}+d_{1}%
+d_{2}+...+d_{i-1}\smallskip$

$\ \ \ \ \ \ \ \ \ \ \ \ \ \ \ \ \ =d_{1}+d_{2}+...+d_{n},$ for all y$_{i}%
$\ $\in$\ Y$_{i}$,

and for 2 $\leq$\ i $\leq$\ n

\ \ \ \ \ sdeg(y$_{i}^{/}$) = $\left|  X_{i}^{/}\right|  =d_{1}+d_{2}%
+...+d_{i-1}\smallskip,$ for all y$_{i}^{/}$\ $\in$\ Y$_{i}^{/}$.

Therefore, signed degree set of G(U, V) is S = $\left\{  d_{1},\sum_{i=1}%
^{2}d_{i},...,\sum_{i=1}^{n}d_{i}\right\}  .$

Clearly, by construction, all the signed bipartite graphs are connected.
Hence, the result follows.\medskip\ \ \ \ \ \ \ \ $\blacksquare$

By interchanging positive edges with negative edges in Theorem 2.1, we obtain
the following result.\medskip

\textbf{Corollary 2.1.} Every set of negative integers is a signed degree set
of some connected signed bipartite graph.\medskip

Finally, we have the following result.\medskip

\textbf{Theorem 2.2.} Every set of integers is a signed degree set of some
connected signed bipartite graph.

\textbf{Proof.} Let S be a set of integers. Then, we have the following five
cases.\newline \textbf{ \ \ \ \ \ \ \ \ \ \ \ }(i)\textbf{ }S is a set of
positive (negative) integers. Then , the result follows by Theorem 2.1
(Corollary 2.1).\newline \textbf{ \ \ \ \ \ \ \ \ \ \ \ }(ii) S = \{0\}. Then,
a signed bipartite graph G(U, V) with $\left|  U\right|  $ = $\left|
V\right|  $ = 2 in which u$_{1}$v$_{1}$, u$_{2}$v$_{2}$ are positive edges and
u$_{1}$v$_{2}$, u$_{2}$v$_{1}$ are negative edges, where u$_{1}$, u$_{2}$
$\in$\ U, v$_{1}$, v$_{2}$ $\in$\ V, has signed degree set S.\newline 

\ \ \ \ \ (iii) S is a set of non-negative (non-positive) integers. Let
$S=S_{1}\cup\{0\},$ where $S_{1}$ is a set of positive (negative) integers.
Then, by Theorem 2.1 (Corollary 2.1), there is a connected signed bipartite
graph $G_{1}(U_{1},V_{1})$ with signed degree set $S_{1}$. Construct a new
signed bipartite graph G(U, V) as follows.

Let $\ \ U=U_{1}\cup\{x_{1}\}\cup\{x_{2}\},\smallskip$

\ \ \ \ \ \ $\ \ V=V_{1}\cup\{y_{1}\}\cup\{y_{2}\},$

with $U_{1}\cap\{x_{i}\}=\phi,\ \{x_{1}\}\cap\{x_{2}\}=\phi,\ V_{1}\cap
\{y_{i}\}=\phi,\ \{y_{1}\}\cap\{y_{2}\}=\phi.$ Let u$_{1}$y$_{1}$, x$_{1}%
$v$_{1}$, x$_{2}$y$_{2}$ be positive edges and u$_{1}$y$_{2}$, x$_{1}$y$_{1}$,
x$_{2}$v$_{1}$ be negative edges, where u$_{1}\in$ U$_{1}$ and v$_{1}\in$
V$_{1}.$ Then, G(U, V) has degree set S. We note that addition of such edges
do not effect the signed degrees of the vertices of G$_{1}$(U$_{1}$, V$_{1}$),
and the vertices x$_{1}$, x$_{2}$, y$_{1}$, y$_{2}$ have signed degrees zero each.\newline 

\ \ \ \ \ (iv) S is a set of non-zero integers. Let $S=S_{1}\cup S_{2}$, where
$S_{1}$ and $S_{2}$ are sets of positive and negative integers respectively.
Then, by Theorem 2.1 and Corollary 2.1, there are connected signed bipartite
graphs G$_{1}$(U$_{1}$, V$_{1}$) and G$_{2}$(U$_{2}$, V$_{2}$) with signed
degree sets $S_{1}$ and $S_{2}$ respectively. Let G$_{1}^{/}$(U$_{1}^{/}$,
V$_{1}^{/}$) and G$_{2}^{/}$(U$_{2}^{/}$, V$_{2}^{/}$) be the copies of
G$_{1}$(U$_{1}$, V$_{1}$) and G$_{2}$(U$_{2}$, V$_{2}$) with signed degree
sets $S_{1}$ and $S_{2}$ respectively. Construct a new signed bipartite graph
G(U, V) as follows.

Let $\ \ U=U_{1}\cup U_{1}^{/}\cup U_{2}\cup U_{2}^{/},\smallskip$

$\ \ \ \ \ \ \ \ V=V_{1}\cup V_{1}^{/}\cup V_{2}\cup V_{2}^{/},$

with $U_{i}\cap U_{j}^{/}=\phi,\ U_{1}\cap U_{2}=\phi,\ U_{1}^{/}\cap
U_{2}^{/}=\phi,\ V_{i}\cap V_{j}^{/}=\phi,\ V_{1}\cap V_{2}=\phi,\ V_{1}%
^{/}\cap V_{2}^{/}=\phi.$ Let u$_{1}$v$_{2}^{/}$, u$_{1}^{/}$v$_{2}$ be
positive edges and u$_{1}$v$_{2}$, u$_{1}^{/}$v$_{2}^{/}$ be negative edges,
where u$_{i}$ $\in$\ U$_{i}$, v$_{i}$ $\in$\ V$_{i}$, u$_{i}^{/}$ $\in
$\ U$_{i}^{/}$ and v$_{i}^{/}$ $\in$\ V$_{i}^{/}$. Then, G(U, V) has signed
degree set S. We note that addition of such edges do not effect the signed
degrees of the vertices of G$_{1}$(U$_{1}$, V$_{1}$), G$_{1}^{/}$(U$_{1}^{/}$,
V$_{1}^{/}$), G$_{2}$(U$_{2}$, V$_{2}$) and G$_{2}^{/}$(U$_{2}^{/}$,
V$_{2}^{/}$).\newline \ \ \ \ \ \ \ \ \ \ \ \ (v) S is a set of all integers .
Let $S=S_{1}\cup S_{2}\cup\{0\}$ , where $S_{1}$ and $S_{2}$ are sets of
positive and negative integers respectively. Then, by Theorem 2.1 and
Corollary 2.1, there exist connected signed bipartite graphs G$_{1}$(U$_{1}$,
V$_{1}$) and G$_{2}$(U$_{2}$, V$_{2}$) with signed degree sets $S_{1}$ and
$S_{2}$ respectively. Construct a new signed bipartite graph G(U, V) as follows.

Let $\ \ U=U_{1}\cup U_{2}\cup\{x\},\smallskip$

\ \ \ \ \ \ $\ \ V=V_{1}\cup V_{2}\cup\{y\},$

with $U_{1}\cap U_{2}=\phi,\ U_{i}\cap\{x\}=\phi,\ V_{1}\cap V_{2}%
=\phi,\ V_{i}\cap\{y\}=\phi.$ Let u$_{1}$v$_{2}$, u$_{2}$y, xv$_{1}$ be
positive edges and u$_{1}$y, u$_{2}$v$_{1}$, xv$_{2}$ be negative edges, where
u$_{i}\in$ U$_{i}$ and v$_{i}\in$ V$_{i}.\ $Then, G(U, V) has signed degree
set S. We note that addition of such edges do not effect the signed degrees of
the vertices of G$_{1}$(U$_{1}$, V$_{1}$) and G$_{2}$(U$_{2}$, V$_{2}$), and
the vertices x and y have signed degrees zero each.\newline 

Clearly, by construction, all the signed bipartite graphs are connected. This
proves the result.\medskip\ \ \ \ \ \ \ \ $\blacksquare\bigskip\bigskip
\bigskip\bigskip$

\begin{center}
\textbf{References\medskip\bigskip}
\end{center}

[1] G. Chartrand, H. Gavlas, F. Harary and M. Schultz, On signed degrees in
signed graphs, Czeck. Math. J., Vol. 44 (1994) 677-690.

[2] S. L. Hakimi, On the realizability of a set of integers as degrees of the
vertices of a graph , SIAM J. Appl. Math. Vol. 10 (1962) 496-506.

[3] F. Harary, On the notion of balance in a signed graph , Michigan Math. J.,
Vol. 2 (1953)143-146.

[4] S. F. Kapoor, A. O. Polimeni and C. E. Wall, Degree sets for graphs, Fund.
Math., Vol. 65 (1977) 189-194.

[5] S. Pirzada and T. A. Naikoo, Signed degree sequences in signed bipartite
graphs, To appear.

[6] S. Pirzada, T. A. Naikoo and F. A. Dar, Degree sets in bipartite and
3-partite graphs, The Oriental J. of Mathematical Sciences, To appear.

[7] S. Pirzada, T. A. Naikoo and F. A. Dar, Signed degree sets in signed
graphs, Czeck. Math. J., To appear.

[8] J. H. Yan , K. W. Lih, D. Kuo and G. J. Chang, Signed degree sequences of
signed graphs, J. Graph Theory, Vol. 26 (1997) 111-117.
\end{document}